\documentclass[a4paper,10pt]{article}
\usepackage{amsmath,amssymb,amsthm}
\usepackage[dvipsnames]{xcolor}
\usepackage{tikz, pgfplots}
\usepackage[unicode, final, colorlinks=true]{hyperref}
\usepackage{textgreek}
\usepackage{enumitem, booktabs}
\usepackage{xcolor}
\usepackage{siunitx}
\usepackage{bbm}
\usepackage{mathtools}
\usepackage{mathabx}
\usepackage{comment}

\hypersetup{citecolor=teal, linkcolor=BrickRed, urlcolor=NavyBlue}

\pgfplotsset{compat=1.17}
\usepgfplotslibrary{groupplots}
\usepgfplotslibrary{fillbetween}
\usetikzlibrary{backgrounds}
\usetikzlibrary{arrows.meta}
\pgfplotsset{plot coordinates/math parser=false}
\usetikzlibrary{external}
\usepgfplotslibrary{patchplots}
\usetikzlibrary{shapes.geometric}
\usetikzlibrary{backgrounds}
\usetikzlibrary{intersections}
\usetikzlibrary{spy}

\newcommand{\R}{\mathbb{R}}
\newcommand{\Z}{\mathbb{Z}}

\newcommand{\ddt}{\partial_t}
\newcommand{\ddx}{\partial_x}
\newcommand{\dds}{\partial_s}

\newcommand{\rhos}{\rho_s}
\newcommand{\cs}{c_s}
\newcommand{\cssqr}{c_s^2}
\newcommand{\csquad}{c_s^4}
\newcommand{\lambdas}{\bar\lambda}

\newcommand{\lambdaf}{\lambda}

\newcommand{\minmod}{\operatorname{minmod}}
\newcommand{\CFL}{\operatorname{CFL}}
\newcommand{\sign}{\operatorname{sign}}

\newcommand{\D}{\mathcal{D}}
\usepackage{bbm}

\newtheorem{theorem}{Theorem}[section]

\newtheorem{remark}[theorem]{Remark}
\newtheorem{algorithm}[theorem]{Algorithm}

\providecommand{\msc}[1]{\textit{2020 MSC:} #1}

\begin{document}
\title{A relaxation approach to the coupling of a two-phase fluid with a linear-elastic solid}
\author{
  Niklas~Kolbe \and
  Siegfried~Müller\footnote{Corresponding author}
}

\date{
  \small
  Institute of Geometry and Applied Mathematics,\\ RWTH Aachen University, Templergraben 55,\\ 52062 Aachen, Germany\\
   \smallskip
   {\tt  \{kolbe,mueller\}@igpm.rwth-aachen.de} \\
   \smallskip
   \today\\
   \bigskip
   \msc{35L65, 35R02, 65M08}
 }
\maketitle

\begin{abstract}
A recently developed coupling strategy for two nonconservative hyperbolic systems is employed to investigate a collapsing vapor bubble embedded in a liquid near a solid. For this purpose, an elastic solid modeled by a linear system of conservation laws is coupled to the two-phase Baer--Nunziato-type model for isothermal fluids, a nonlinear hyperbolic system with non-conservative products. For the coupling of the two systems the Jin-Xin relaxation concept is employed and embedded in a second order finite volume scheme.
For a proof of concept simulations in one space dimension are performed.
\end{abstract}

\section{Introduction}

The dynamics of cavitation bubbles and their interaction with a compliant wall are of interest
in numerous real world applications arising, for instance, in biology, engineering,
and medical applications, such as
cavitation erosion of under water structures
\cite{Philipp-Lauterborn:98},
lithotripsy and sonoporation
\cite{Ohl-Wolfrum:03,Johnson-Colonius-Cleveland:09},
and
cavitation-enhanced ablation of materials, e.g., in biological tissues
\cite{Brujan-Nahen-Schmidt-Vogel:01b}.
In particular, the transient distribution of 
stresses is important in order to understand 
the cause of
the observed damage comparing breaking points of the material.

To numerically investigate  
cavitation erosion, a model of 
two-phase compressible fluids can be coupled with a model 
for elastic solids.
It is well-known from continuum mechanics that at the interface the negative stress tensor of the material model equals the pressure of the fluid, also the displacement and  velocity in normal direction coincide, see~\cite[Chap.~5.1]{BazilevsTakizawaTezduyar-2013}. 
These transition conditions couple the dynamics of the fluid and the material. 

In the context of fluid-structure interaction (FSI) typically iterative approaches 
are applied to couple the models until the coupling condition is satisfied up to a prescribed tolerance before proceeding to the next time step.
For a review on the topic we refer to \cite{HouWangLayton:12}.
To avoid such iterative techniques and the alternating application of two solvers strategies have been developed that aim to solve a coupled Riemann problem at the interface.
fluid-structure interaction (FSI), 
For instance, Wang et al.~\cite{Wang-Rallu-Gerbeau-Farhat:11}
solve a half-Riemann problem in the fluid at the interface. Here the fluid state at the interface is determined 
by moving along the rarefaction curve
and shock curve, respectively, starting from an adjacent fluid state until the normal velocity coincides with the normal interface velocity. 
From this fluid state the flow induced load is determined by means of the transmission condition for the normal stress component.
Alternatively, 
Banks et al.~\cite{Banks-Henshaw-Sjogreen:13} 
developed a coupling strategy 
based on the solution of a one-dimensional coupled Riemann problem for two hyperbolic systems of conservation laws at the interface, 
where both the fluid equations and the linear elastic equations have the \emph{same} number of unknowns. The solution of the coupled Riemann problem provides
ghost states to compute the numerical fluxes at the interface for the fluid and the solid solver, respectively.

In the hyperbolic community coupling strategies for two hyperbolic fluid systems have been developed and analyzed extensively for many years, see for instance 
\cite{AmbrosoChalonsCoquelGodlewskiLagoutiereRaviartSeguin-2008,ChalonsRaviartSeguin-2008,
GodlewskiRaviart2004aa,GodlewskiRaviart2005aa, GodlewskiLe-ThanhRaviart2005aa} and references therein. 
The basic idea is to interpret the coupled problem as two boundary value problems. The boundary values on either side of the coupling interface are then chosen such that some coupling condition is satisfied, e.g., continuity of some state variables or fluxes at the coupling interface, which is referred to either as state coupling or flux coupling, respectively. The main task is then to find admissible boundary conditions at the coupling interface. For this purpose, a coupled Riemann problem is solved with initial data taken from both sides of the interface. Finally, in the numerical discretization, e.g., by the finite volume or the discontinuous Galerkin method, for each system a flux is computed at the coupling interface from the boundary values solving the coupled Riemann problem. 
Similar ideas are also investigated in the context of coupling of hyperbolic transport dynamics on networks.
Exemplarily, we refer to results by Bressan et al.~\cite{BressanCanicGaravello2014aa}.

As the aforementioned approach is specific to the coupling of fluids an alternative approach has been developed in \cite{HertyMuellerGerhardXiangWang:2018,Sikstel:794674}. The main idea is to solve two coupled half-Riemann problems at the interface, where the coupling states on either side are spanned by the Lax curves of the corresponding system. The positions on the Lax-curves are then determined by plugging the coupling states into the coupling conditions.
This concept has been modified in \cite{herty2023centr,herty2023centrschemtwo} in a way that avoids the computation of the Lax curves for nonlinear systems. To determine the coupling states at the interface, the Jin-Xin relaxation procedure \cite{JinXin1995aa} is applied separately to each of the systems. The resulting linear relaxation systems are then coupled by means of two linear half-Riemann problems proceeding as before. Performing the relaxation limit the coupling states for the original problem are determined. 
Recently, this strategy has been extended to the coupling of two nonconservative hyperbolic systems, see \cite{kolbe2023numerical}. 

In the following we apply this strategy to the coupling of a linear-elastic solid with a two-phase fluid. To this end we summarize in Sect.~\ref{sec:problem} the coupling problem consisting of the two models for the solid and the fluid as well as the coupling conditions. To determine appropriate coupling states at the interface we then consider the coupling of the corresponding Jin-Xin relaxation systems in Sect.~\ref{sec:relaxation}. The numerical scheme based on a finite volume discretization is introduced in Sect.~\ref{sec:num}. Eventually, numerical simulations for a one-dimensional vapor bubble embedded in a liquid collapsing near a solid are presented in Sect.~\ref{sec:experiments}.

\section{Coupling two-phase flow to a linear-elastic structure}\label{sec:problem}
This section introduces the coupling problem that we address in this work. We present a model for a linear-elastic structure,
a Baer--Nunziato type model of two-phase flow 
and formulate coupling conditions.
\paragraph{Model of a linear-elastic solid.}
We adopt the one-dimensional elastic structure model
\begin{subequations}\label{eq:elastic}
\begin{align}
  \ddt w - \frac{1}{\rho_s} \ddx \sigma &= 0 && (t,x) \in \R^+ \times \R^-,\\
  \ddt \sigma - \rho_s \cssqr \ddx w &= 0  && (t,x) \in \R^+ \times \R^-,
\end{align}
\end{subequations}
  for the deformation velocity $w$ and the shear stress $\sigma$. Material parameters are the density of the material $\rho_s$ and the dilatation wave velocity $\cs$.
  For more details on the derivation of linear-elastic structure models we refer to \cite{dickopp:hal-01121991, HertyMuellerGerhardXiangWang:2018}.

\paragraph{Two-phase flow model.}
To model two-phase flow we employ the 
barotropic Baer--Nunziato model for two phases.
\begin{subequations}\label{eq:BN}
  In this model the volume fraction of the first phase, $\alpha_1$, is governed by
\begin{equation}
  \ddt \alpha_1 + v_I \ddx \alpha_1 = S^v_{\alpha,k} + S^p_{\alpha,k} \qquad (t,x) \in \R^+ \times \R^+,
\end{equation}
where $v_I$ denotes the interfacial velocity, and the one of the second phase, $\alpha_2$, by $\alpha_1 + \alpha_2=1$.
Conservation of mass and momentum give rise to
\begin{align}
\ddt (\alpha_k \rho_k) + \ddx (\alpha_k \rho_k v_k) &= S^v_{\alpha \rho,k} + S^p_{\alpha\rho,k} && (t,x) \in \R^+   \times \R^+, \\
  \ddt (\alpha_k \rho_k v_k) + \ddx (\alpha_k \rho_k v_k^2 + \alpha_k p_k) + p_I \ddx \alpha_\ell &= S^v_{\alpha \rho v,k} + S^p_{\alpha\rho v,k} && (t,x) \in \R^+ \times \R^+,
\end{align}
\end{subequations}
which hold for $k\in\{1,2\}$. By $\rho_k$, $v_k$ and $p_k=p_k(\rho_k)$ we refer to the density, the velocity and the pressure of phase $k$, respectively. The index of the opposite phase $\ell=\ell(k)$ is taken $2$ if $k=1$ and $1$ if $k=2$ and $p_I$ denotes the interfacial pressure. We define the mixture density, momentum and pressure by
\begin{equation}
\rho = \alpha_1 \rho_1 + \alpha_2 \rho_2, \quad \rho v = \alpha_1 \rho_1 v_1 + \alpha_2 \rho_2 v_2, \quad p = \alpha_1 p_1 + \alpha_2 p_2.
\end{equation}
The interfacial velocity and pressure can be taken e.g., as the linear combinations
\begin{equation}
\label{eq:interfacial-states}
v_I = \beta_1 v_1 + \beta_2 v_2, \qquad p_I = \beta_2 p_1 + \beta_1 p_2
\end{equation}
where
\[
  \beta_k = \frac{d_k \alpha_k \rho_k}{d_1 \alpha_1 \rho_1 + d_2 \alpha_2 \rho_2}, \qquad k \in \{ 1, 2\}
  \]
  for chosen constant parameters $d_1, d_2 \in [0, 1]$ such that $d_1+d_2=1$. Note that choosing $d_k=1$ and $d_\ell=0$ we obtain $v_I = v_k$ and $p_I = p_\ell$. We assume here that both phases correspond to an isothermal stiffened gas with equation of state given by
  \begin{equation}
  \label{eq:eos}
    p_k(\rho_k) = c_k^2 \rho_k -\pi_k,\qquad k \in \{1, 2\},
  \end{equation}
  where $c_k$ and $\pi_k$ denote the isothermal speed of sound and the minimal pressure corresponding to phase $k$, respectively.
  
  The model accounts for velocity and pressure relaxation by means of the source terms on the right-hand sides of \eqref{eq:BN} determined by 
  \begin{subequations}
  \label{eq:relaxation}
  \begin{align}
    & S^v_{\alpha ,k} =0  ,\ S^v_{\alpha \rho,k} = 0,\ S^v_{\alpha \rho v,k}=\theta_v \alpha_k \rho_k (v-v_k),\\
    & S^p_{\alpha ,k} =\theta_p \alpha_k  (p-p_k),\ S^p_{\alpha \rho,k} = 0,\ S^p_{\alpha \rho v,k}=0
  \end{align}
  \end{subequations}
  with inverse velocity relaxation time $\theta_v$ and inverse pressure relaxation time $\theta_p$.
  For more details on the Baer-Nunziato-type model for isothermal fluids we refer to~\cite{HantkeMuellerSikstelThein:2024-arxiv}.
  
  \paragraph{Coupling condition.}
At the interface we impose the equality of all velocities:
\begin{subequations}\label{eq:cplvelocity}
\begin{align}
  w(0^-,t) &= v_1(0^+, t),  \label{eq:cplvelocitya}\\
  w(0^-, t) &= v_2(0^+, t) \label{eq:cplvelocityb}
\end{align}
\end{subequations}
both, for a.e.~$t\geq 0$. In addition, the following relation between solid stress and fluid pressures should be satisfied for a.e.~$t\geq 0$:
\begin{subequations}\label{eq:cplstress}
\begin{align}
    \sigma(0^-,t) &= -p_1(0^+, t),  \label{eq:cplstressa}\\
  \sigma(0^-, t) &= -p_2(0^+, t). \label{eq:cplstressb}
\end{align}
\end{subequations}

\section{The coupled Jin-Xin relaxation system}\label{sec:relaxation}
In this section we consider a Jin-Xin type relaxation of the coupling problem from Section~\ref{sec:problem}.

\subsection{Jin-Xin relaxation of the linear-elastic solid model}
We first consider the Jin-Xin relaxation of \eqref{eq:elastic} and introduce the vector notation
\begin{align}
\bar U = (w, \sigma)^T \in \bar{ \mathcal{D}} \subset \R^2,  \qquad \bar F(\bar U) = 
- \left( \frac{\sigma}{\rho_s}, \rho_s \cssqr w \right)^T
\end{align}
as well as the new variable
\[
\bar V = (V^w, V^\sigma)^T \in \R^2.
\]
Then, for any relaxation rate $\varepsilon>0$ the Jin-Xin relaxation system of \eqref{eq:elastic} takes the form
\begin{equation}\label{eq:relelastic}
\begin{aligned}
  \ddt \bar U + \ddx \bar V &= 0 && (t,x) \in \R^+ \times \R^-, \\
  \ddt \bar V + \bar \lambda^2  \, \ddx \bar U &= \frac{1}{\varepsilon} \left(\bar F(\bar U) - \bar V \right)
                                             && (t,x) \in \R^+ \times \R^-.
\end{aligned}
\end{equation}
For stability of the relaxation the scalar $\bar \lambda>0$ must be taken such that the subcharacteristic condition
is satisfied, see \cite{chen1994hyperconserlaws}.

\subsection{Jin-Xin relaxation of the two phase flow model}
As \eqref{eq:BN} is a nonconservative hyperbolic system 
when neglecting the relaxation terms on the right-hand side we employ the recent generalized approach from \cite{kolbe2023numerical} for the Jin-Xin relaxation. Therefore, we define
\begin{align*}
  U &= (\alpha_1, \alpha_1 \rho_1, \alpha_1 \rho_1 v_1, \alpha_2 \rho_2, \alpha_2 \rho_2 v_2)^T \in \mathcal{D} \subset \R^5, \\
  F(U) &= \left(0 , \alpha_1 \rho_1 v_1,  \alpha_1 \rho_1 v_1^2 +  \alpha_1 p_1,  \alpha_2 \rho_2 v_2,  \alpha_2 \rho_2 v_2^2 +  \alpha_2 p_2 \right)^T,\\
  g(U) &= (v_I, 0, p_I, 0, -p_I)^T, \qquad G(U) = \begin{pmatrix} g(U) & 0 & 0& 0 &0 \end{pmatrix}.
\end{align*}
Since it holds $\alpha_2 = 1 - \alpha_1$, we can now write the homogeneous version of~\eqref{eq:BN} as
\[
  \ddt U + \ddx F(U) + G(U) \ddx U = 0
\]
or equivalently in quasilinear form
\[
  \ddt U + A(U) \ddx U = 0
\]
for $A(U)= DF(U) + G(U)$. Since the solution concept of nonconservative systems is ambiguous at discontinuities a class of weak solutions needs to be fixed. To this end we follow \cite{dal1995defin} and consider a family of paths
\[
\Phi : [0,1] \times \D \times \D \to \D,
\]
which, along with certain continuity properties, see e.g., \cite{pares2006numer}, satisfies $\Phi(0, U^-, U^+) = U^-$ and $\Phi(1, U^-, U^+) = U^+$. A simple example is the family of segment paths given by $\Phi(s, U^-, U^+) = U^- + s (U^+ - U^-)$. Next, we define the path integral
\begin{equation}
  \mathcal{P}(U^-, U^+; A) = \int_0^1 A(\Phi(s; U^-, U^+) \dds \Phi(s; U^-, U^+) \, ds
\end{equation}
and the distribution
\begin{equation}\label{eq:pathdistribution}
  \begin{split}
    \langle [A(U(t, \cdot)) \, \ddx U(t, \cdot)]_\Phi, \phi \rangle  &= \int_0^\infty A(U(t, x))\,\ddx U(t, x) \phi(x) dx \\
      &\quad + \sum_{x_D \in D(t)} \mathcal{P}(U(t,x_D^-), U(t, x_D^+); A) \phi(x_D)
  \end{split}
\end{equation}
for piecewise smooth functions $U(t, \cdot)$ with discontinuities at time $t$ contained in the set $D(t)\subset \R^+$ and denoting by $U(t,x_D^-)$ and $U(t, x_D^+)$ the limit of $U(t, \cdot)$ as the discontinuity $x_D$ is approached from the left and the right, respectively. We note that \eqref{eq:pathdistribution} gives rise to a class of distributional solutions of \eqref{eq:BN}.

Moreover, the distribution \eqref{eq:pathdistribution} is used in our relaxation approach, for which we define the nonlocal operator
\begin{equation}\label{eq:T}
  T[U(t, \cdot)](x) \coloneqq - \langle [A(U(t, \cdot)) \, \ddx U(t, \cdot)]_{\Phi}, \mathbbm{1}_{[x, \infty)} \rangle.
\end{equation}
This operator can be roughly interpreted as integral over the nonconservative product $A(U) \ddx U$ from $\infty$ to $x$. For brevity, we neglect the dependence on $t$ and $x$ in the notation in the following. Assuming the existence of a state $U_\infty$ such that $U-U_\infty$ has compact support in $\R^+$ the operator \eqref{eq:T} can be rewritten as
\[
  T[U] = F(U) - F(U_\infty) - \langle [G(U(t, \cdot)) \, \ddx U(t, \cdot)]_{\Phi}, \mathbbm{1}_{[x, \infty)} \rangle.
  \]
  For the relaxation we introduce the new variable
  \[
    V \coloneqq (V^{\alpha_1}, V^{\alpha_1 \rho_1}, V^{\alpha_1 \rho_1 v_1}, V^{\alpha_2 \rho_2},V^{\alpha_2 \rho_2 v_2})^T \in \R^5.
  \]
  Then we state the following Jin-Xin relaxation system for \eqref{eq:BN}:
\begin{equation}\label{eq:relBN}
  \begin{aligned}
    \ddt U + \ddx V &= 0 && (t,x) \in \R^+ \times \R^+,\\
    \ddt V + \lambda^2 \, \ddx  U &= \frac{1}{\varepsilon} \left(T[U] - V \right) && (t,x) \in \R^+ \times \R^+.
  \end{aligned}
\end{equation}
For simplicity, we consider a simultaneous relaxation of both systems, \eqref{eq:elastic} and \eqref{eq:BN}, in which we employ the same relaxation rate $\varepsilon$ in \eqref{eq:relelastic} and \eqref{eq:relBN}. For stability the subcharacteristic condition with respect to both systems needs to be satisfied, which is given, if we take
\[
  \bar \lambda = \max \{ \hat \rho(D \bar F(\bar U)):~\bar U\in \bar \D \}, \qquad  \lambda = \max \{\hat \rho( A( U)):~ U\in \D \}
  \]
  with $\hat \rho$ denoting the spectral radius.

  \subsection{Coupling conditions for the relaxed problem}\label{sec:relcoupling}
  In addition to \eqref{eq:cplvelocity} and \eqref{eq:cplstress}, further conditions are required for the well-posedness of the coupling problem given by the relaxed systems \eqref{eq:relelastic} and \eqref{eq:relBN}. As we aim for a solution of the limit problem from Section~\ref{sec:problem} consistency with respect to the original conditions in the relaxation limit is required, see~\cite{herty2023centrschemtwo}. Apart from the needed well-posedness of the resulting Riemann solver and the consistency the choice of additional conditions is free; they can be chosen such that the coupling problem does not become too complex and allows for a computable solution.

  \begin{figure}
    \centering
\includegraphics{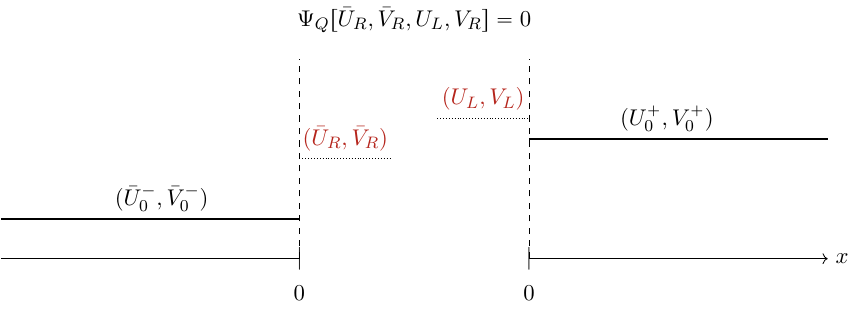}
\caption{The coupled relaxation system at the interface. The coupling states $\bar U_R$, $\bar V_R$, $U_L$ and $V_L$ constitute suitable boundary data with respect to the (numerical) traces $\bar U_0^-$, $\bar V_0^-$, $U_0^+$, $V_0^+$ and satisfy a coupling condition encoded in $\Phi_Q$.}\label{fig:couplingstates}
\end{figure}

Due to the nonlocality in \eqref{eq:relBN} we construct conditions for the general case that there is a jump between the upper trace at the coupling interface that we denote by $U_0$ and $V_0$ and the imposed left boundary states $U_L$ and $V_L$. Moreover, we denote the right boundary states, for which the conditions are imposed for the negative half axis by $\bar U_R$ and $\bar V_R$. In the following we also use analogue notations for the corresponding components. In Figure~\ref{fig:couplingstates} we visualize the relation between trace and coupling states.

  Focussing first on the relaxed elastic system \eqref{eq:relelastic} we note that consistency in the relaxation limit requires $\bar V_R = \bar F(\bar U_R)$. Analyzing the components of this equation allows us to rewrite the coupling states as
  \begin{equation}\label{eq:relelasticconsistency}
    \sigma_R = - \rho_s V^w_R , \qquad 
    w_R = - \frac{V^\sigma_R}{\cssqr \rho_s}.
  \end{equation}

  For the relaxed Baer--Nunziato system on the positive half-axis we follow \cite{kolbe2023numerical} to get the consistency condition
  \begin{equation}\label{eq:relBNconsistency}
    \begin{split}
      V_R &= V_0 - \mathcal{P}(U_L, U_0^+; A)\\
      &= V_0^+ + F(U_L) - F(U_0^+) -  \mathcal{P}(U_L, U_0^+; G)
    \end{split}
  \end{equation}
  Considering the second and fourth components of \eqref{eq:relBNconsistency} allows us to rewrite the velocities of the phases as
  \begin{equation}\label{eq:relBNvconsistency}
    v_{kL} = \frac{V^{\alpha_k \rho_k}_L - V^{\alpha_k \rho_k}_0 + \alpha_{k0} \rho_{k0} v_{k0}}{\alpha_{kL} \rho_{kL}} \qquad k \in \{1,2\}
  \end{equation}
  noting that the corresponding components of $\mathcal{P}(U_L, U_0^+; G)$ are zero.

  As the interfacial pressure depends on $\rho_1$ and $\rho_2$, we use the notation $p_I = p_I(U)$. The fifth component of $\mathcal{P}(U_L, U_0^+; G)$ can then be written as
  \[
    \mathcal{P}_{p_I} (U_L, U_0^+) \coloneqq  \int_0^1 p_I(\Phi(s, U_L, U_0^+)) \dds \phi_1(s, U_L ,U_0^+) \, ds
  \]
  with $\phi_1$ denoting the first component of $\Phi$. From components three and five of \eqref{eq:relBNconsistency} we now obtain
  \begin{subequations}\label{eq:relBNpconsistency}
  \begin{align}
    \alpha_{1L} p_{1L} &= V^{\alpha_1 \rho_1 v_1}_L - V^{\alpha_1 \rho_1 v_1}_0 + \alpha_{10} \rho_{10} v_{10}^2 + \alpha_{10}p_1(\rho_{10}) - \mathcal{P}_{p_I} (U_L, U_0^+) - \alpha_{1L} \rho_{1L} v_{1L}^2,\\
    \alpha_{2L} p_{2L} &= V^{\alpha_2 \rho_2 v_2}_L - V^{\alpha_2 \rho_2 v_2}_0 + \alpha_{20} \rho_{20} v_{20}^2 + \alpha_{20}p_2(\rho_{20}) + \mathcal{P}_{p_I} (U_L, U_0^+) - \alpha_{2L} \rho_{2L} v_{2L}^2.
  \end{align}
\end{subequations}

To derive two additional coupling conditions we combine the original condition $w_R = v_{1L} = v_{2L}$ with \eqref{eq:relelasticconsistency} and \eqref{eq:relBNvconsistency}. In this way we obtain
\begin{equation}\label{eq:relcouplingv}
    - \frac{V^\sigma_R}{\cssqr \rho_s} =
    \frac{V^{\alpha_k \rho_k}_L - V^{\alpha_k \rho_k}_0 + \alpha_{k0} \rho_{k0} v_{k0}}{\alpha_{kL} \rho_{kL}} \qquad k \in \{1,2\}.
  \end{equation}

  The coupling conditions \eqref{eq:cplstress} imply, in particular, that $-\sigma_R = \alpha_1 p_1 + \alpha_2 p_2$. Substituting \eqref{eq:relelasticconsistency} and \eqref{eq:relBNpconsistency} into this equation  we obtain another coupling condition, which after exploiting consistency with respect to $V_R^w$ reads
  \begin{equation}\label{eq:relcouplingp}
    \sigma_R =  \sum_{k=1}^2 \left( V^{\alpha_k \rho_k v_k}_L - V^{\alpha_k \rho_k v_k}_0 + \alpha_{k0} \rho_{k0} v_{k0}^2 + \alpha_{k0}p_k(\rho_{k0}) - \alpha_{kL} \rho_{kL} v_{kL}^2\right).
  \end{equation}
  By construction the new coupling conditions become redundant with the original ones from Section~\ref{sec:problem} in the relaxation limit $\varepsilon \to 0$.
  \begin{remark}
    For a.e.~$t\geq0$ the solution of the continuous coupling problem has no jumps at the interface, i.e., we have $U_0^+=U_L$ and $V_0^+=V_L$ and from \eqref{eq:relcouplingv} and \eqref{eq:relcouplingp} we get
    \begin{subequations}
      \begin{align}
        V^\sigma(0^-,t) &= -\cssqr \rho_s v_1(0^+, t) = -\cssqr \rho_s v_2(0^+, t),\\
        \rho_s V^w(0^-, t) &= \alpha_1 p_1 (\rho_1(0^-, t)) + \alpha_2 p_2 (\rho_2(0^-, t))
      \end{align}
    \end{subequations}
    for a.e.~$t\geq 0$.
  \end{remark}
  \begin{remark}
We note that the path integrals do not enter the coupling conditions \eqref{eq:cplvelocity}, \eqref{eq:cplstress}, \eqref{eq:relcouplingv} and \eqref{eq:relcouplingp} of the coupled relaxation system. Since from a physical point of view the interface is a material interface where a solid and a fluid are coupled, this is in agreement with the findings in \cite{dal1995defin,LeFloch:1988}. There it is proven that in nonconservative hyperbolic systems the generalized Rankine-Hugoniot jump conditions depend on the path in case of shock waves, i.e., a genuinely nonlinear characteristic field, but are independent of the path in case of contact waves, i.e., a linearly degenerated characteristic field.
\end{remark}

\section{Numerical scheme}\label{sec:num}
In this section we introduce a finite volume scheme for the coupled two-phase fluid structure problem stated in Section \ref{sec:problem}. We first elaborate on the general discretization in Section \ref{sec:fv} and focus on the numerical solution of the coupling problem in Section \ref{sec:RS}.

\subsection{Finite volume scheme}\label{sec:fv}
The Lax--Friedrichs like scheme that we employ is based on the approach in \cite{herty2023centr, kolbe2023numerical}.  It has been derived from an upwind discretization of the Jin-Xin relaxation systems \eqref{eq:relelastic} and \eqref{eq:relBN} in the relaxation limit $\varepsilon \to 0$.

\begin{figure}
  \includegraphics{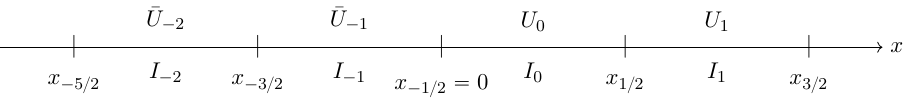}
  \caption{Discretization of the real line as used in our finite volume scheme.}\label{fig:linedisc}
\end{figure}
We partition the real line into the uniform mesh cells $I_j=(x_{j-1/2}, x_{j+1/2})$ of width $\Delta x$ so that the origin is located at the cell interface $x_{-1/2}$ as shown in Figure~\ref{fig:linedisc}. Additionally, we introduce the time instances $t^n= n \Delta t$ for the uniform time step size satisfying the condition
\begin{equation}
\label{eq:time-discretization}
  \Delta t= \CFL \frac{\Delta x}{\max \{ \bar \lambda, \lambda \}}
\end{equation}
given a suitable Courant number $0<\CFL \leq 1$. Next, we consider volume averages over the cell $I_j$ at time $t^n$ for the state variables in both, system \eqref{eq:elastic} and system \eqref{eq:BN} denoted using the corresponding indices and combine them into the discrete states
\begin{equation}
W_j^n = \bar U_j^n = \begin{pmatrix} w_j^n \\ \sigma_j^n \end{pmatrix} \quad \text{for }j<0, \qquad
W_j^n =  U_j^n = \begin{pmatrix} \alpha_{1,j}^n \\ \alpha_{1,j}^n \, \rho_{1,j}^n  \\ \alpha_{1,j}^n \, \rho_{1,j}^n \, v_{1,j}^k \\
  \alpha_{2,j}^n \, \rho_{2,j}^n  \\ \alpha_{2,j}^n \, \rho_{2,j}^n \, v_{2,j}^k\end{pmatrix} \quad \text{for }j \geq 0.
\end{equation}
The vector $W_j^n$ here is introduced for notational convenience and depending on the index $j$ refers to either the state of the elastic structure or the multiphase-fluid. Our scheme takes the form
\begin{equation}\label{eq:scheme}
W_j^{n+1} = W_j^n - \frac{\Delta t}{ \Delta x } \left(\mathcal{F}_{j+1/2}^{n,-} - \mathcal{F}_{j+1/2}^{n,+} \right), \qquad j \in \Z
\end{equation}
with numerical fluxes depending on the spatial index. Away from the interface the scheme is conservative and we have
\begin{equation}
  \begin{aligned}
    \mathcal{F}_{j+1/2}^{n, \text{struc}} \coloneqq \mathcal{F}_{j+1/2}^{n,-} = \mathcal{F}_{j+1/2}^{n,+} = \frac 12 \left( \bar F(\bar U_{j-1}^n)  + \bar F(\bar U_{j}^n) \right) &- \frac{\bar \lambda}{2} (\bar U_j^n - \bar U_{j-1}^n) \\
    &- \frac{\Delta x}{2} (\bar S_j^{n,-} - \bar S_{j-1}^{n,+})
  \end{aligned}, \quad j<0
\end{equation}
for the elastic structure left and
\begin{equation}\label{eq:numfluxfluid}
\begin{aligned}
    \mathcal{F}_{j+1/2}^{n, \text{fluid}} \coloneqq \mathcal{F}_{j+1/2}^{n,-} = \mathcal{F}_{j+1/2}^{n,+} = \frac 12 \left( T_{j-1}[ U^n]  +  T_j [ U^n] \right) &- \frac{\lambda}{2} ( U_j^n -  U_{j-1}^n) \\
    &- \frac{\Delta x}{2} (S_j^{n,-} - S_{j-1}^{n,+})
  \end{aligned}, \quad j>0
\end{equation}
for the multiphase fluid right from the interface. In the numerical flux~\eqref{eq:numfluxfluid} we employ a discetized form of operator \eqref{eq:T} given by
\begin{equation}
T_j[U^n] = F(U_j^n) - F(U_\infty) - \sum_{k \geq j} \mathcal{P}(U_k^n, U_{k+1}^n; G).
\end{equation}
We note that for $j>0$ scheme \eqref{eq:scheme} can be written as a path conservative scheme employing only local contributions to update the state in a given cell, see \cite{kolbe2023numerical}.

For high resolution we use the MUSCL scheme \cite{vanleerUltimateConservativeDifference1977} making use of slope reconstructions of the characteristic variables. To prevent oscillatory behavior we use the minmod limiter, which has led to accurate simulations under moderate time step size restrictions. The reconstructed slopes take the form
\begin{align}
  \bar S_j^{n, \pm} &\coloneqq \minmod\left( \frac{\bar F(\bar U_j^n) - \bar F(\bar U_{j-1}^n) \pm \bar \lambda (\bar U_j^n - \bar U_{j-1}^n)  }{2 \Delta x},
  \frac{\bar F(\bar U_{j+1}^n) - \bar F(\bar U_{j}^n) \pm \bar \lambda (\bar U_{j+1}^n - \bar U_{j}^n)  }{2 \Delta x} \right),\\
    S_j^{n, \pm} &\coloneqq \minmod\left( \frac{T_j[U^n] - T_{j-1}[U^n] \pm  \lambda (U_j^n - U_{j-1}^n)  }{2 \Delta x},
                   \frac{T_{j+1}[U^n] -  T_j[U^n] \pm \lambda (U_{j+1}^n - U_{j}^n)  }{2 \Delta x} \right),
\end{align}
where the \emph{minmod} operator is given by
\[
  \minmod(a,b) =
  \begin{cases*}
    0 & if $\sign(a) \neq \sign(b)$ \\
    a & if $|a| \leq |b|$ and $\sign(a) = \sign(b)$ \\
    b & if $|a| > |b|$ and $\sign(a) = \sign(b)$
  \end{cases*}
\]
for scalar arguments $a,b \in \R$ and component-wise in case of vectors.

At the interface we employ the numerical flux
\begin{equation}\label{eq:leftinterfaceflux}
  \mathcal{F}_{-1/2}^{n,-} = \frac 12 \left( \bar F(\bar U_{-1}^n)  + \bar V_R^n \right) - \frac{\bar \lambda}{2} (\bar U_R^n - \bar U_{-1}^n)
\end{equation}
from the negative half-axis and
\begin{equation}\label{eq:rightinterfaceflux}
  \mathcal{F}_{-1/2}^{n,+} = \frac 12 \left( V_L^n  + T_0[U^n] \right) - \frac{\lambda}{2} (U_0^n - U_{L}^n)
\end{equation}
from the positive half-axis.

Note that we do not account for velocity and pressure relaxation in the above finite volume discretization. Instead, we perform an operator splitting in the fluid where after each update of the fluid we perform an instantanuous velocity and pressure relaxation to equilibrium, i.e., the source terms of \eqref{eq:BN} are considered in the limit as both, $\theta_v\to \infty$ and $\theta_p \to \infty$. For this purpose, we apply Algorithms 5.1 and 5.2 in \cite{HantkeMuellerSikstelThein:2024-arxiv} by which the state in each cell is modified. Thus, after each time step the velocities and pressures are in local equilibrium, i.e., $v_1=v_2=v_I=v^\infty$ and $p_1=p_2=p_I=p^\infty$.

\subsection{Riemann solver}\label{sec:RS}
This section is concerned with the computation of the quantities $\bar U_R^n$, $\bar V_R^n$, $U_L^n$ and $V_L^n$ emerging in the numerical interface fluxes \eqref{eq:leftinterfaceflux} and \eqref{eq:rightinterfaceflux}. These so-called coupling states will be determined using a Riemann solver (RS) for the coupled Jin-Xin relaxation system discussed in Section~\ref{sec:relaxation}, which is a map of the form
\begin{equation}\label{eq:RS}
  \mathcal{RS}: (\bar U_{-1}^n, \bar V_{-1}^n, U_{0}^n, V_0^n) \mapsto (\bar U_R^n, \bar V_R^n, U_L^n, V_L^n).
\end{equation}
Since the coupling conditions for the coupled relaxation have been designed with consistency to the original problem in mind this RS can be used within the scheme \eqref{eq:scheme}. To determine the missing coupling states the input variables to \eqref{eq:RS} with respect to the auxiliary states $\bar V$ and $V$ are replaced with the corresponding relaxation limits, i.e., we take
\begin{equation}
   \bar V_{-1}^n = \bar F( \bar U_{-1}^n) \qquad \text{and} \qquad V_0^n = T_0[U^n].
 \end{equation}
 As only a fixed time instance is concerned we neglect the index $n$ in the remainder of this section.

 The RS~\eqref{eq:RS} provides a solution of the half-Riemann problem in the sense that a) the coupling states are suitable boundary data for the hyperbolic problem on both sides of the coupling and b) the coupling data satisfy the coupling conditions that we derived in Section~\ref{sec:relcoupling}. Point a) is addressed by connecting the coupling states to the trace states, i.e., the arguments of \eqref{eq:RS} by waves with velocities having the correct sign, see \cite{garavello2006traffflownetwor}. Since the relaxation system is linear and can be decoupled component-wise this leads to the relations
 \begin{equation}\label{eq:laxrelation}
   \bar V_R = \bar V_{-1} + \bar \lambda (\bar U_{-1} - \bar U_R) \qquad \text{and} \qquad  V_L = V_{0} +  \lambda ( U_L - U_{0}),
 \end{equation}
 see~\cite{kolbe2023numerical} for details.

 Regarding b) we first employ~\eqref{eq:laxrelation} within the coupling conditions \eqref{eq:cplvelocity}, \eqref{eq:cplstress}, \eqref{eq:relcouplingv} and \eqref{eq:relcouplingp} to obtain a nonlinear system in $\bar U_R$ and $U_L$, which thus has 7 scalar unknowns. To determine a solution of this system we first rewrite coupling condition \eqref{eq:cplstressa} as
 \begin{equation}\label{eq:rssigma}
   \sigma_R = \frac{\pi_1 \alpha_{1L} - c_1^2 [\alpha_1 \rho_1]_{L}}{\alpha_{1L}}
 \end{equation}
 and \eqref{eq:cplvelocitya} as
 \begin{equation}\label{eq:rsmomentum1}
[\alpha_1 \rho_1 v_1]_L = w_R [\alpha_1 \rho_1]_L.
 \end{equation}
 Substituting \eqref{eq:rssigma} in \eqref{eq:relcouplingv} for $k=1$ allows us to express $\alpha_{1L}$ as a function in $[\alpha_1\rho_1]_{L}$, i.e., $\alpha_{1L}([\alpha_1 \rho_1]_{L})$. Similarly, employing \eqref{eq:cplstressb} and substituting \eqref{eq:rssigma} and $\alpha_{1L}([\alpha_1 \rho_1]_{L})$ we obtain a formula for $[\alpha_2 \rho_2]_{L}$ in $[\alpha_1 \rho_1]_{L}$, i.e., $[\alpha_2 \rho_2]_{L}([\alpha_1\rho_1]_{L})$. Combining the new expressions with coupling condition \eqref{eq:relcouplingv} for $k=2$ we obtain $R_1([\alpha_1 \rho_1]_{L})=0$ for a third order polynomial 
 $R_1$.

Finally, we propose the following 7 step procedure to define the RS \eqref{eq:RS} and thus determine the coupling states:
\begin{algorithm}[Computation of coupling states]\label{algo:coupling}
~\\[-4mm]
 \begin{enumerate}
 \item Compute the roots of $R_1(x)$. Out of the real roots $x_1<x_2<x_3$ we select $[\alpha_1 \rho_1]_{L}=x_2$.
 \item Determine $[\alpha_2 \rho_2]_{L}([\alpha_1\rho_1]_{L})$.
 \item Determine $\alpha_{1L}([\alpha_1 \rho_1]_{L})$.
 \item Determine $\sigma_R$ inserting $\alpha_{1L}$ and $[\alpha_1 \rho_1]_{L}$ into \eqref{eq:rssigma}.
 \item Substitute $[\alpha_1 \rho_1]_{L}$ and $\sigma_R$ into \eqref{eq:relcouplingv} to obtain the quadratic equation $R_2([\alpha_2 \rho_2 v_2]_{L})=0$. The smaller root of $R_2$ defines $[\alpha_2 \rho_2 v_2]_{L}$.
 \item Determine $w_R$ inserting $[\alpha_2 \rho_2 v_2]_{L}$ and $[\alpha_2 \rho_2]_{L}$ into \eqref{eq:cplvelocityb}.
 \item Determine $[\alpha_1 \rho_1 v_1]_L$ inserting $[\alpha_1 \rho_1]_{L}$ and $w_R$ into \eqref{eq:rsmomentum1}.
 \end{enumerate}
 \end{algorithm}
 The detailed formulas for $\alpha_{1L}([\alpha_1 \rho_1]_{L})$, $[\alpha_2 \rho_2]_{L}([\alpha_1\rho_1]_{L})$, $R_1$ and $R_2$ are provided in Appendix \ref{sec:algoappendix}.
 \begin{remark}
   In all numerical tests with application relevant parameters and trace states the polynomial $R_1$ has always had three real roots. Moreover, the smallest root $x_1$ has been consistently negative and thus $[\alpha_1 \rho_1]_L = x_1$ would have been unphysical. Similarly, choosing $[\alpha_1 \rho_1]_L = x_3$ has consistently led to an unphysical volume fraction $\alpha_{1L} \notin [0,1]$ in step 3.
Furthermore, in our tests the larger root of $R_2$ has been several magnitudes larger than the corresponding trace state, i.e.~choosing the larger root has led to $[\alpha_2 \rho_2 v_2]_{L} \gg [\alpha_2 \rho_2 v_2]_{0}$. This has motivated our choice in steps 1 and 5 of the procedure within the RS.
 \end{remark}

 \section{Numerical experiments}\label{sec:experiments}
 In this section we apply our method to numerically investigate the collapse of a one-dimensional vapor bubble embedded in a liquid near a linear-elastic solid. This configuration is motivated by the two-dimensional simulations of a collapsing hot gas bubble in cold water performed in~\cite{Sikstel:794674}.

\paragraph{Configuration.}
The spatial domain $\Omega=(-0.2,0.2)$ m is split into the two parts $\Omega_s=(-0.2,0)$ m and $\Omega_f=(0,0.2)$ m. To the left of the material interface $x=0$ we assume a linear-elastic material with zero dilatation velocity ($w=0$ m/s) and shear stress $\sigma=-3.5\times 10^7$~N/m$^2$. Here we consider steel characterized by the parameters $\rho_s=7800$ kg/m$^3$ and $\cs=5990$ m/s.
Right of the material interface we consider an isothermal two-phase fluid
with an almost pure vapor state inside the bubble placed in $\Omega_B = (0.075, 0.125)$ m determined by
$\alpha_1=0.9$, $p_1=p_2= 3.5 \times 10^3$ N/m$^2$, $v_1=v_2=0$ m/s 
embedded in an almost pure water state characterized by
$\alpha_1=0.1$, $p_1=p_2=5\times 3.5 \times 10^6$ N/m$^2$, $v_1=v_2=0$ m/s within $\Omega_f \setminus \Omega_B$. 
The phases $k=1$ and $k=2$ correspond to water vapor and liquid water respectively, and are modeled by an ideal isothermal gas law and an isothermal stiffened gas law with parameters 
$c_1=367.58$ m/s, $\pi_1=0$ N/m$^2$ and $c_2=1483.3$ m/s and $\pi_2=1.1358 \times 10^9$ N/m$^2$ in \eqref{eq:eos} corresponding to a constant reference temperature $T=T_{ref}=293.15$ K.
The interfacial velocity and pressure are taken as $v_I=v$ and $p_I=p$, respectively, motivated by the local mechanical equilibrium that we assume.
Note that both states are at rest and in local velocity and pressure equilibrium. At the material interface the coupling conditions \eqref{eq:cplvelocity} and \eqref{eq:cplstress} are satisfied.

\paragraph{Discretization.}
The solid and the fluid part of the computational domain are uniformly discretized using $\Delta x = 1/3000$ m, and thus $600$ cells in both, $\Omega_s$ and $\Omega_f$. For the solid and fluid part we use the same time discretization $\Delta t$ using the CFL number $\CFL=0.2$ in \eqref{eq:time-discretization}. On the domain boundaries at $x=-0.2$m and $x=0.2$ we impose homogeneous Neumann boundary conditions.

 \begin{figure}
   \centering
   \includegraphics{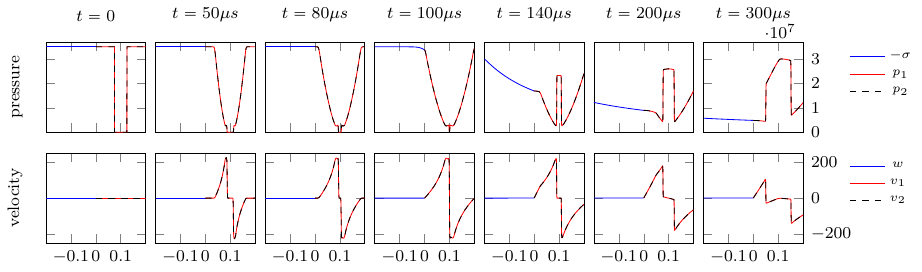}
    \caption{Numerical solution in terms of pressure and velocity showing seven time instances until $t =300 \mu s$ over the computational domain with coupling interface at $x=0$, where the structure and the fluids interact.}\label{fig:pv1}
  \end{figure}

   \begin{figure}
     \centering
     \includegraphics{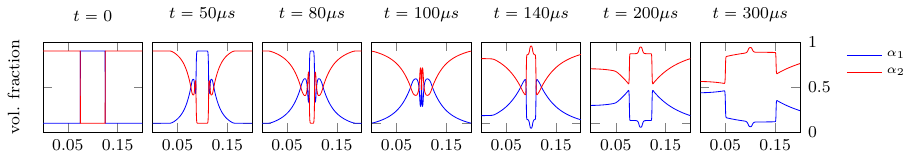}
    \caption{Volume fractions of water vapor ($\alpha_1$) and liquid water ($\alpha_2$) in the numerical solution until $t =300 \mu s$ over the multiphase fluid domain $[0, 0.2]$.}
    \label{fig:vfractions1}
  \end{figure}

   \begin{figure}
     \centering
     \includegraphics{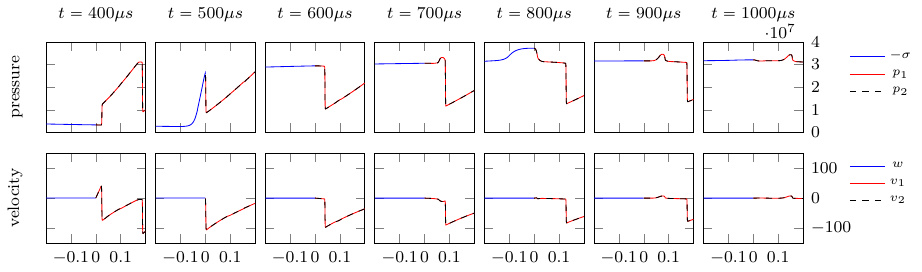}
    \caption{Numerical solution in terms of pressure and velocity showing seven time instances from $t= 400 \mu s$ until $t =1000 \mu s$ over the computational domain with coupling interface at $x=0$, where the structure and the fluids interact.}\label{fig:pv2}
  \end{figure}

     \begin{figure}
       \centering
       \includegraphics{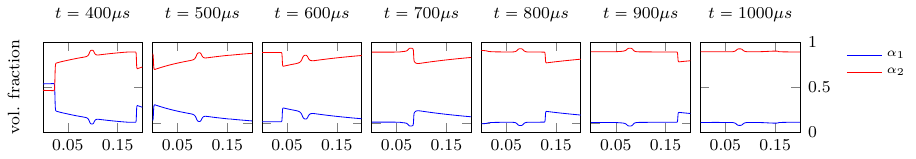}
    \caption{Volume fractions of water vapor ($\alpha_1$) and liquid water ($\alpha_2$) in the numerical solution from $t =400 \mu s$ until $t=1000 \mu s$ over the multiphase fluid domain $[0, 0.2]$.}\label{fig:vfractions2}
  \end{figure}

  \paragraph{Results.}
In Figures \ref{fig:pv1} and \ref{fig:pv2} we present the phasic fluid velocities and pressures as well as the deformation velocity and the negative shear stress in the solid, respectively, for different time instances. We note that the fluid velocities and pressures are in local equilibrium due to the instantanuous relaxation. Furthermore, there are no layers at the coupling interface (material interface) confirming the jump conditions  \eqref{eq:cplvelocity} and \eqref{eq:cplstress}.

In Figures \ref{fig:vfractions1} and \ref{fig:vfractions2} the corresponding volume fractions in the fluid are presented. These change because of the instantanuous pressure relaxation. We emphasize that there is no change in the phasic mass $\alpha_1\rho_1$ and $\alpha_2\rho_2$ because we do not account for phase change due to relaxation of Gibbs free energies.

The simulation can be split into several parts. Initially at time $t=0$, two-phase Riemann problems are solved at the liquid-vapor interfaces of the vapor bubble consisting of a rarefaction wave and a shock wave. The rarefaction waves run into the liquid water whereas the shock waves run in opposite directions into the vapor bubble followed by the liquid-vapor interface, see $t=50$ $\mu$s and $t=80$ $\mu$s. At time $t=100$ $\mu$s the leading edge of the left-running rarefaction wave interacts with the solid-fluid interface at $x=0$ where it is being partially transmitted into the solid and partially reflected into the fluid due to the coupling conditions. Note that the transmitted wave is running much faster in the solid than in the reflected wave in the fluid because of the significantly larger characteristic speed in the solid, see time $t=140$ $\mu$s.
In the meantime the two shock waves as well as  the liquid-vapor material waves  are focussing in the center of the bubble where they are  reflected, see times $t=140,\ 200,\ 300,\ 400$ $\mu$s. 
At time $t=500$ $\mu$s the left running reflected shock and liquid-vapor material wave interact with the solid at the coupling interface. There the waves are partially reflected into the fluid and partially transmitted into the the solid, see time $t=600,\ 700$ $\mu$s. 
At time $t=700$ $\mu$s the reflected liquid-vapor material wave interacts with the remainder of the original bubble. This causes a wave in the pressure that is interacting with the liquid-solid material interface, see time $t=800$ $\mu$s where it is partially transmitted and reflected, see times $t=800,\ 900,\ 1000$~$\mu$s.

\begin{figure}
  \centering
  \includegraphics{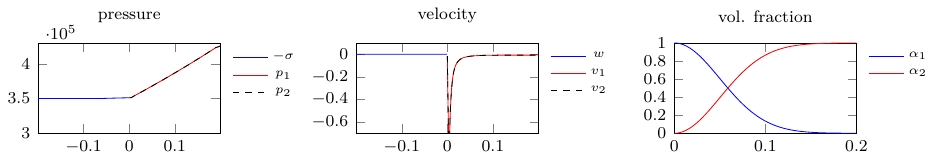}
    \caption{Numerical solution of the grid convergence experiment in terms of pressure, velocity and volume fraction on $2 \times 600$ mesh cells at the time instance $t=10 \mu s$.}\label{fig:gridstudyini}
  \end{figure}

\paragraph{Grid convergence.} To validate our approach we investigate the convergence of the numerical solution as the grid is refined in a simplified numerical experiment. 
We employ the same material parameters and computational domain as above but the modified constant shear stress $\sigma=-3.5\times 10^5$~N/m$^2$ for the linear-elastic material on $\Omega_s$ as well as the volume fraction $\alpha_1= e^{-200 x^2}$ and the fluid pressures $p_1=p_2 = 3.5\times 10^5 e^x $N/m$^2$ for the two-phase fluid model on $\Omega_f$. These modifications have been made to ensure continuous initial data. In our analysis we compare numerical solutions at the time $t=10 \mu$ s at which an decrease of shear stress and fluid velocity is visible at the interface, see Figure~\ref{fig:gridstudyini}.

\begin{table}
  \caption{Coupling errors \eqref{eq:cplerrors} with respect to velocity and pressure and EoCs under grid refinement.}\label{tab:cpl} \vspace{5pt}
  \centering
  \begin{tabular}{rlrlr}
    \toprule
$N$ & $E_C^1$ & EoC & $E_C^2$ & EoC\\ \midrule
200 & 4.461e-02 &  & 9.584e+01 &  \\
400 & 2.506e-02 & 0.831 & 4.865e+01 & 0.978\\
800 & 1.323e-02 & 0.922 & 2.461e+01 & 0.983\\
1600 & 6.514e-03 & 1.022 & 1.191e+01 & 1.047\\ \bottomrule
\end{tabular}
\end{table}

Using our approach we compute solutions on various grids; to reduce the numerical error with respect to time we thereby choose the time steps according to $\Delta t = \CFL \Delta x^2/\bar \lambda$.  In Table~\ref{tab:cpl} we show the discrete coupling errors
\begin{equation}\label{eq:cplerrors}
E_C^1 \coloneqq \left|\frac{[\alpha_1 \rho_1 v_1]_0 + [\alpha_2 \rho_2 v_2]_0}{[\alpha_1 \rho_1]_0 + [\alpha_2 \rho_2]_0} - w_{-1}\right| \quad \text{and}\quad  E_C^2 \coloneqq |[\alpha_1]_{0} p_1([\rho_1]_0) + [\alpha_2]_{0} p_2([\rho_2]_0) + \sigma_0 | 
\end{equation}
in terms of velocity and pressure at time $t$. These errors are derived from \eqref{eq:cplvelocity} and \eqref{eq:cplstress}. Along we show the experimental order of convergence (EoC\footnote{The EoC is computed by the formula $\text{EoC} = \log_2(E_1/E_2)$ with $E_1$ and $E_2$ denoting the error in two consecutive lines of the table.}). The coupling errors decrease as the grid is refined and the EoCs indicate a first order convergence with respect to the grid resolution.

\begin{table}
  \caption{Discrete $L^1$ errors with respect to dilatation velocity, shear stress, mixture density and momentum under grid refinement.\\}\label{tab:L1}
  \centering
\begin{tabular}{rlrlrlrlr} \toprule
$N$ & $E_w$ & EoC & $E_\sigma$ & EoC & $E_\rho$ & EoC & $E_{\rho v}$ & EoC\\ \midrule
200 & 3.372e-07 &  & 1.521e+01 &  & 7.599e-03 &  & 8.304e-03 & \\
400 & 1.617e-07 & 1.060 & 7.416e+00 & 1.036 & 2.043e-03 & 1.895 & 3.651e-03 & 1.186\\
800 & 6.148e-08 & 1.395 & 2.867e+00 & 1.371 & 5.472e-04 & 1.900 & 1.534e-03 & 1.250\\
1600 & 1.956e-08 & 1.652 & 9.160e-01 & 1.646 & 1.474e-04 & 1.893 & 5.748e-04 & 1.417\\ \bottomrule
\end{tabular}
\end{table}

In addition we present in Table~\ref{tab:L1} the discrete $L^1$ errors for various model states such as dilatation velocity, shear stress, mixture density and momentum over the computational domain at the final time. With an exception of the mixture momentum the corresponding EoCs tend to two as the grid is refined thanks to the second order MUSCL reconstructions we employ. Note that due to the nonsmooth profile shown in Figure~\ref{fig:gridstudyini} EoCs closer to 2 are not expected in this experiment.

\appendix
\section{Solution of the nonlinear system}\label{sec:algoappendix}
Here we provide the detailed formulas used in Algorithm~\ref{algo:coupling} to obtain the coupling states and solve the corresponding nonlinear system.

Writing the polynomial $R_1$ in step 1 as $R_1(x) = a_3 x^3 + a_2 x^2 + a_1 x + a_0$, its coefficients are given by
 \begin{align*}
          & a_3 =  ((\lambdaf -w_0) \rhos (([\alpha_1 \rho_1]_0 c_1^2 +[\alpha_2 \rho_2]_0 c_2^2) \lambdaf -[\alpha_1 \rho_1 v_1]_0 c_1^2 -[\alpha_2 \rho_2 v_2]_0 c_2^2) \cssqr - \\
          & \qquad \lambdas ((c_1^2 (\pi_2 -\sigma_0) [\alpha_1 \rho_1]_0 +[\alpha_2 \rho_2]_0 c_2^2 (\pi_1 -\sigma_0)) \lambdaf -c_1^2 (\pi_2 -\sigma_0) [\alpha_1 \rho_1 v_1]_0 -[\alpha_2 \rho_2 v_2]_0 c_2^2 (\pi_1 -\sigma_0))) \lambdas,\\
          & a_2 = (\rhos^2 (\lambdaf -w_0)^{2} \csquad-\lambdas \rhos (([\alpha_1 \rho_1]_0 c_1^2 +[\alpha_2 \rho_2]_0 c_2^2 -2 \sigma_0 +\pi_1 +\pi_2) \lambdaf -[\alpha_1 \rho_1 v_1]_0 c_1^2 + \\
          &\qquad (-\pi_1 -\pi_2 +2 \sigma_0) w_0 -[\alpha_2 \rho_2 v_2]_0 c_2^2) \cssqr +\lambdas^2 (\pi_2 -\sigma_0) (\pi_1 -\sigma_0)) (\lambdaf [\alpha_1 \rho_1]_0 -[\alpha_1 \rho_1 v_1]_0),\\
          & a_1 = \texttt{-}2 \rhos \cssqr (\lambdaf [\alpha_1 \rho_1]_0 -[\alpha_1 \rho_1 v_1]_0)^{2} (\rhos (\lambdaf -w_0) \cssqr -\frac{\lambdas (\pi_1 +\pi_2 -2 \sigma_0)}{2}),\\
          & a_0 =  \rhos^2 \csquad (\lambdaf [\alpha_1 \rho_1]_0-[\alpha_1 \rho_1 v_1]_0)^{3} .
        \end{align*}

        In step 2 we express $[\alpha_2 \rho_2]_L$ as a function of $[\alpha_1 \rho_1]_L$ via
         \begin{align*}
      &[\alpha_2\rho_2]_L([\alpha_1\rho_1]_L) = -\frac{N_1([\alpha_1\rho_1]_L)}{D_1([\alpha_1\rho_1]_L)},\\
      &N_1([\alpha_1\rho_1]_L) = (\lambdaf [\alpha_1\rho_1]_L \cssqr \rhos -\lambdaf [\alpha_1 \rho_1]_0 \cssqr \rhos -[\alpha_1\rho_1]_L w_0 \cssqr \rhos +\lambdas [\alpha_1\rho_1]_L \sigma_0 -\pi_2 \lambdas [\alpha_1\rho_1]_L + \\
      &\qquad [\alpha_1 \rho_1 v_1]_0 \cssqr \rhos) (\lambdaf [\alpha_1\rho_1]_L \cssqr \rhos -\lambdaf [\alpha_1 \rho_1]_0 \cssqr \rhos +\lambdas [\alpha_1\rho_1]_L^2 c_1^2 -[\alpha_1\rho_1]_L w_0 \cssqr \rhos +\lambdas [\alpha_1\rho_1]_L \sigma_0 \\
      & \qquad -\lambdas [\alpha_1\rho_1]_L \pi_1 +[\alpha_1 \rho_1 v_1]_0 \cssqr \rhos),\\
      &D_1([\alpha_1\rho_1]_L) = [\alpha_1\rho_1]_L \lambdas c_2^2 (\lambdaf [\alpha_1\rho_1]_L \cssqr \rhos -\lambdaf [\alpha_1 \rho_1]_0 \cssqr \rhos -[\alpha_1\rho_1]_L w_0 \cssqr \rhos + \\
      & \qquad\lambdas [\alpha_1\rho_1]_L \sigma_0 -\lambdas [\alpha_1\rho_1]_L \pi_1 +[\alpha_1 \rho_1 v_1]_0 \cssqr \rhos).
         \end{align*}
         
The first volume fraction at the interface in step 3 can be expressed as
      \begin{align*}
        &\alpha_{1L}([\alpha_1\rho_1]_L) = -\frac{N_2([\alpha_1\rho_1]_L)}{D_2([\alpha_1\rho_1]_L)},\\
        &N_2([\alpha_1\rho_1]_L) = \lambdas [\alpha_1\rho_1]_L^2 c_1^2,\\
        &D_2([\alpha_1\rho_1]_L) = \lambdaf [\alpha_1\rho_1]_L \cssqr \rhos -\lambdaf \,[\alpha_1 \rho_1]_0 \cssqr \,\rhos -[\alpha_1\rho_1]_L \,w_0 \cssqr \,\rhos + \\
        &\qquad \lambdas \,[\alpha_1\rho_1]_L \,\sigma_0 -\lambdas \,[\alpha_1\rho_1]_L \,\pi_1 +[\alpha_1 \rho_1 v_1]_0 \cssqr \,\rhos.        
      \end{align*}

      Writing the quadratic polynomial in step 5 as $R_2(x) = b_2 x^2 +  b_1 x + b_0 $ its coefficients are
            \begin{align*}
        & b_2 = -([\alpha_1\rho_1]_L +[\alpha_2\rho_2]_L \,) [\alpha_2 \rho_2]_0  [\alpha_1 \rho_1]_0 ,\\
        & b_1 = (\lambdas \,\rhos +\lambdaf  [\alpha_2\rho_2]_L +\lambdaf \,[\alpha_1\rho_1]_L) [\alpha_1 \rho_1]_0  [\alpha_2\rho_2]_L [\alpha_2 \rho_2]_0,\\
        & b_0 = (((-\pi_1 +\pi_2) \alpha_{10} -([\alpha_1 \rho_1 v_1]_0 +[\alpha_2 \rho_2 v_2]_0) \lambdaf -w_0 \,\rhos \,\lambdas +\sigma_0 -\pi_2) [\alpha_2\rho_2]_L  [\alpha_2 \rho_2]_0 + \\
        & \qquad \lambdas \,\rhos \,[\alpha_2 \rho_2]_0 +[\alpha_2 \rho_2 v_2]_0^2 [\alpha_2\rho_2]_L +[\alpha_2\rho_2]_L \,c_2^2 \,[\alpha_2 \rho_2]_0^2-[\alpha_2 \rho_2]_0 \,\lambdas \,\rhos) [\alpha_1 \rho_1]_0 [\alpha_2\rho_2]_L + \\
         & \qquad([\alpha_1 \rho_1]_0^2 c_1^2 +[\alpha_1 \rho_1 v_1]_0^2) [\alpha_2\rho_2]_L^2 [\alpha_2 \rho_2]_0.
            \end{align*}
            \paragraph{Funding.} The authors thank the Deutsche Forschungsgemeinschaft (DFG, German Research Foundation) for the financial support under 320021702/GRK2326 (Graduate College Energy, Entropy, and Dissipative Dynamics), through SPP 2410 (Hyperbolic Balance Laws in Fluid Mechanics: Complexity, Scales, Randomness) within the Project 525842915, through SPP 2311 (Robust Coupling of Continuum-Biomechanical In Silico Models to Establish Active Biological System Models for Later Use in Clinical Applications – Co-Design of Modelling, Numerics and Usability) within the Project 548864771 and though Project 461365406 (New traffic models considering complex geometries and data). Support by the ERS Open Seed Fund of RWTH Aachen University through project OPSF781 is also acknowledged.
            
            {\scriptsize
              \bibliographystyle{abbrvurl}
              \bibliography{references.bib}
            }
            
\end{document}